\documentclass[11pt]{article}
\usepackage{amsmath,amsthm,amsfonts,latexsym,amssymb,enumerate,color}
\usepackage{graphicx}
\usepackage{mathrsfs}

\def\GA{\mathcal A}
\def\CC{\mathbb C}
\def\GC{\mathcal C}
\def\DD{\mathbb D}
\def\GD{\mathcal D}

\def\GF{\mathcal F}
\def\G{\mathcal G}
\def\GH{\mathcal H}
\def\GJ{\mathcal J}

\def\GL{\mathcal L}
\def\GM{\mathcal M}
\def\GN{\mathcal N}

\def\GR{\mathcal R}
\def\RR{\mathbb R}

\def\TT{\mathbb T}
\def\GX{\mathcal X}

\def\ker{\mathop{\rm ker}\nolimits}

\def\rank{\mathop{\rm rank}\nolimits}

\def\q{\quad}

\newtheorem{thm}{Theorem}[section]
\newtheorem{prop}[thm]{Proposition}
\newtheorem{defn}[thm]{Definition}
\newtheorem{lem}[thm]{Lemma}
\newtheorem{cor}[thm]{Corollary}
\newtheorem{rem}[thm]{Remark}

\numberwithin{equation}{section}

\def\beginpf{\begin{proof}}
\def\endpf{\end{proof}}
\def\beq{\begin{equation}}
\def\eeq{\end{equation}}

\def\ol{\overline}
\def\til{\tilde}
\def\wt{\widetilde}

\def\sigab{\Sigma_{a,b}}

\def\H2b{\overline{H^2_0}}

\begin{document}

\title{Kernels of paired operators and their adjoints}

\author{M.~Cristina C\^amara,\thanks{
Center for Mathematical Analysis, Geometry and Dynamical Systems,
Instituto Superior T\'ecnico, Universidade de Lisboa, 
Av. Rovisco Pais, 1049-001 Lisboa, Portugal.
 \tt ccamara@math.ist.utl.pt} 
 \and  Jonathan R.~Partington\thanks{School of Mathematics, University of Leeds, Leeds LS2~9JT, U.K. {\tt j.r.partington@leeds.ac.uk}
}
\  }

\date{}

\maketitle

\begin{abstract}
We review the basic properties of paired operators and their adjoints, the transposed paired operators, with particular reference to commutation relations, and we study the properties of their kernels, bringing out their similarities and also, somewhat surprisingly, their stark differences. Various notions expressing different invariance properties are also reviewed and we extend to paired operators some known invariance results.
\end{abstract}



\let\thefootnote\relax\footnote{Work partially funded by FCT/Portugal through project UIDB/04459/2020 with DOI identifier 10-54499/UIDP/04459/2020.}


\section{Introduction}

The canonical example of a singular integral equation on $L^p (\TT)$, $1<p<\infty$, is
\beq\label{1.1}
a(t)f(t)+b(t)(S_\TT f)(t)=g(t)
\eeq
where $a$ and $b$ are bounded functions in $\TT$, $g$ is a given function in $L^p (\TT)$ and $S_\TT$ is the singular integral operator with Cauchy kernel defined by
\beq\label{1.2}
(S_\TT f)(t)=\frac{1}{\pi i} \hbox{PV} \int_\TT{\frac{f(\tau)}{\tau - t} d\tau}.
\eeq

The boundedness, Fredholmness and invertibility of the operator
\beq\label{1.3}
T=a I + b S_\TT
\eeq
defined by the left hand side of \eqref {1.1}, as well as methods to obtain explicit solutions to that equation and systems of equations of the same type, have been extensively studied, mainly in connection with factorisation of Wiener--Hopf type (see for example \cite{widom, GK, MP} and references therein). In general certain conditions are imposed on $a$ and $b$, such as $a^2 - b^2$ being bounded away from zero \cite{widom, MP}. For $p=2$, it was Shinbrot who first derived a constructive general method to find explicit solutions to \eqref{1.1},   by placing it in an abstract Hilbert space setting, in \cite{shinbrot}, using the fact that
\beq\label{1.4}
I = P^+ + P^- \qquad \hbox{and} \qquad S_\TT = P^+-P^-,
\eeq
where $P^\pm$ are the orthogonal projections from $L^2$ onto $H_+^2 := H^2(\DD)$ and $H_-^2 = (H_+^2)^\perp$, respectively. Indeed \eqref{1.4} implies that \eqref{1.3} can be written as a {\em paired operator} \cite{MP}
\beq\label{1.5}
S_{A,B} = AP + BQ
\eeq
where $A$ and $B$ are bounded operators and $P$, $Q$ are complementary projections on a Banach space. In \cite{shinbrot} Shinbrot assumed that $A^{-1}$ existed and $A^{-1}B$ was selfadjoint and of single sign, noting that the situation would be much more complicated in the general case.

Naturally, one can also consider singular integral equations of the form
\beq\label{1.6}
af + b S_\TT f = g
\eeq
to which, by \eqref{1.4}, one can associate in a similar way another paired operator of the form
\beq\label{1.7}
\Sigma_{A,B} = P^+ A + P^- B,
\eeq
which we call a {\em transposed paired operator}.

If $A$ and $B$ are invertible operators, then \eqref{1.5} and \eqref{1.7} are said to be of 
 normal type \cite{MP}. This is by far the most studied case. Then we have that any operator of the form \eqref{1.7} is equivalent to an operator of the form \eqref{1.5}, i.e., they are related by products with invertible operators:
\beq\label{1.8}
PA + QB = [(I+P A B^{-1} Q) A] (B^{-1} P + A^{-1} Q) [(I - Q A B^{-1} P) B].
\eeq
Moreover, if $AB = BA$, then \cite[III, Theorem 1.2]{MP} 
\beq\label{1.9}
PA + QB = [(I + P A B^{-1} Q) B^{-1}] (AP + BQ) [(I - Q A B^{-1} P) B].
\eeq
This is true, in particular, if $A$ and $B$ are multiplication operators on $L^2$, as it happens in the case of singular integral operators.

Thus, for many purposes (such as the study of boundedness, Fredholmness, invertibility, solvability of equations) it suffices to consider only one type of paired operators. This may explain why transposed paired operators have hardly been considered separately in the literature on this topic, in spite of renewed interest in recent years  \cite{gu,speck,NY10, NY14, DDS24, CGP23, CP24}.

In this paper we shall work in the context $X = L^2(\TT)$, where $\TT$ is the unit circle, and $P = P^+$, $Q = P^-$. For functions $a, b \in L^\infty(\TT) =: L^\infty$, we take $A = M_a$ and $B = M_b$, the corresponding multiplication operators. Then we write $S_{a,b}$ and $\Sigma_{a,b}$ for the corresponding paired and transposed  paired operators, i.e.,
\beq\label{1.10}
S_{a,b}= aP^+ + b P^- \qquad \hbox{and} \qquad \Sigma_{a,b}=P^+ aI + P^- bI.
\eeq

Writing $T_a$ for the Toeplitz operator $P^+ a P^+$ on $H^2_+$, we see that
\[
P^+ S_{a,b}{}_{| H^2_+} = P^+ \Sigma_{a,b}{}_{|H^2_+} = T_a,
\]
so that paired operators are dilations of Toeplitz operators. Moreover, if $b \in G L^\infty$ (that is, invertible in $L^\infty$), then $S_{a,b}$ and $\Sigma_{a,b}$ are equivalent, and both are equivalent after extension to $T_{a/b}$ \cite{Cam17} (see also \cite{BTsk} for more on these concepts). Hence we have an isomorphism of kernels:
\beq\label{1.11}
\hbox{ker} \Sigma_{a,b} \approx \hbox{ker} S_{a,b} \approx \hbox{ker} T_{a/b}.
\eeq

The papers \cite{CGP23} and \cite{CP24} were mainly concerned with paired operators $S_{a,b}$, but  the theory of transposed paired operators is also very rich, and we shall also be discussing these, bringing out some of the similarities and differences.

We shall assume throughout the paper  that $a,b \in L^\infty$ and
\beq\label{1.12}
a \ne 0, \quad b \ne 0 \quad \hbox{a.e. on } \TT, \quad \hbox{and} \quad a = b \quad \hbox{or} \quad a-b \ne 0 \quad \hbox{a.e. on }\TT.
\eeq
If $a \ne b$ then we say that the pair $(a,b)$ is {\em nondegenerate}.

After describing some of the basic properties of paired operators in Section~\ref{sec:2} we shall present some of the properties of paired kernels in Section~\ref{sec:3}, before moving on to transposed paired kernels in Section~\ref{sec:4}, in which the similarities and differences between the two notions will be identified. Finally, Section~\ref{sec:5} is devoted to reviewing various notions expressing different invariance properties, such as almost-invariant and nearly-invariant subspaces for a given operator, analysing their relations and extending some known invariance results to kernels of paired operators. Along the way we get new insights into Toeplitz and Hankel operators and new formulations of known results for Toeplitz kernels.

Some notation: we shall write $\tilde T_a$ for the dual Toeplitz  operator $P^- a P^-$ on $H_-^2$ and $\tilde H_a$ for the Hankel operator $P^+ a P^-: H_-^2 \rightarrow H_+^2$. In addition, we shall say that $O_- \in H_-^2$ (resp. $\overline{H^\infty}$) is {\em outer} if $\overline{z} \overline{O_-}$ (resp. $\overline{O}_-$) is {\em outer} in $H_+^2$ (resp. $H^\infty$).

\section{Basic properties}
\label{sec:2}

We start by noting that the two types of paired operator, $S_{a,b}$ and $\Sigma_{a,b}$,
are closely connected. Indeed, if $b \in \G L^\infty$, then $
S_{a,b}$ and $\Sigma_{a,b}$ are equivalent, i.e., they
can be related by invertible operators 
\beq
\Sigma_{a,b} = E S_{a,b} F
\eeq
with $E,F$ invertible \cite[Ch.~III, Sec.~1]{MP}.

Moreover, for any $a,b \in L^\infty$ we have
\beq\label{eq:16jul2.2}
S_{a,b}^*=\Sigma_{\ol a,\ol b}
\eeq
and
\beq
\ol{S_{a,b}f}= z S_{\ol b,\ol a} \ol z \ol f = z \Sigma^*_{b,a}(\ol z \ol f),
\eeq
 so that some properties of
$\Sigma_{a,b}$ can be deduced easily from the corresponding properties
of $S_{a,b}$. 
For example, we have the following optimal estimates from \cite{CGP23}.

\begin{thm}
For $a,b \in L^\infty$, let $m=\max(\|a\|_\infty,\|b\|_\infty)$
and $M=\|a\|_\infty + \|b\|_\infty$. Then 
\[
m \le \|S_{a,b}\|=\|\Sigma_{\ol a,\ol b}\| \le \min(M, \sqrt 2\, m).
\]
\end{thm}
Clearly, the same inequalities hold for $\Sigma_{a,b}$.
As an immediate consequence we have the following result, but we can also give
an easy direct proof.
\begin{cor}\label{cor:14aug2.1A}
\beq
S_{a,b}=0 \iff \Sigma_{a,b}=0 \iff a=b=0,
\eeq
and it follows that the symbol pair $(a,b)$ is unique for any operator of the form 
\eqref{1.10}. 
\end{cor}

\beginpf
If $S_{a,b}=0$, then, applying this to functions $\phi_+ \in H^2_+$ we see that $a\phi_+=0$ for all
such functions and  so $a=0$. Similarly, applying this to functions in $H^2_-$ we see
that $b=0$. By \eqref{eq:16jul2.2} we obtain the same conclusion for $\Sigma_{a,b}$.
The converse is obvious.
\endpf

Since, for any $\eta \in L^\infty$,
\[
M_\eta = S_{\eta,\eta}=\Sigma_{\eta,\eta},
\]
we have that $S_{a,b}$ (or $\Sigma_{a,b}$) is a multiplication operator on $L^2$ if and only if $a=b$.

In the same way that nonzero Toeplitz operators cannot be compact (let alone finite rank),
a similar result holds for paired operators, and hence for transposed paired operators.

\begin{prop}
The operators $S_{a,b}$ and $\Sigma_{a,b}$
can only be compact if $a=b=0$.
\end{prop}
\beginpf
Note that $\Sigma_{a,b}=P^+aI+P^-bI$, and this can
only be compact if   the operators   $P^+aI$ and $P^-bI$
are both compact, because a Hilbert space sequence $(x_n+y_n)$ with $\langle x_n,y_m\rangle=0$
for all $n$ and $m$ converges if and only if the two sequence $(x_n)$ and $(y_n)$ both
converge. 

Now $P^+aI: L^2 \to H^2_+$ is an extension of the Toeplitz operator
$T_a$ on $H^2_+$ and so is not compact unless $a=0$. Similarly for $P^-bI$.
Hence $\Sigma_{a,b}$ is not compact unless $a=b=0$, and so neither is $S_{a,b}=\Sigma_{\ol a,\ol b}^*$.
\endpf

We shall require certain commutation results from \cite[Prop.~3.1]{CGP23}.
Since for two operators on a Hilbert space we have $XT=TX$ if and only
if $T^*X^* = X^* T^*$, and indeed
the adjoint of the commutator $[X,T]=XT-TX$ is the commutator
$[T^*,X^*]$, which has the same rank, we can immediately apply the results
we quote to the context of transposed paired kernels.
Similarly, questions on the commutant of two Toeplitz
operators, as analysed   by Brown and Halmos in \cite{BH63},
have their generalizations to paired operators.
The first part of the following result is due to Gu \cite[Thm.~3.1]{gu}
and the second part follows by duality (this is a slight improvement on  
 \cite[Lem.~3.1]{CGP23}).

\begin{prop}\label{prop:16jul2.3}
For any nondegenerate $(a,b)$, $(\til a,\til b)$ satisfying \eqref{1.12}, we have:\\
(i) $S_{a,b}S_{\til a,\til b}= S_{c,d} \iff \tilde a \in H^\infty \hbox{ and } \tilde b \in \ol{H^\infty}$, and in that case $c=a\til a$ and $d=b\til  b$;\\
(ii) $\Sigma_{a,b}\Sigma_{\til a,\til b}=\Sigma_{c,d} \iff a \in \ol{H^\infty}
 \hbox{ and }   b \in H^\infty$,
 and in that case $c=a\til a$ and $d=b \til b$.
 \end{prop}
 
 As an immediate consequence of Proposition~\ref{prop:16jul2.3} we see the following.
 
 \begin{cor}
 For nondegenerate $(a,b)$ and $(\til a,\til b)$ we have
 $S_{a,b}S_{\til a,\til b}=0$ if and only if $\til a \in H^\infty$, $\til b \in \ol{H^\infty}$
 and $a\til a=b \til b=0$.
 \end{cor}
 
More generally, we have by direct calculation (see also   \cite[Prop.~3.2]{CGP23}):
 
 \begin{prop}\label{prop:2.2jun14}
 (i) $S_{a,b}S_{\til a,\til b}-S_{a\til a,b\til b}=
 (a-b)(P^+\til b P^--P^- \til a P^+)$;\\
 (ii) $\Sigma_{a,b}\Sigma_{\til a,\til b}-\Sigma_{a\til a,b\til b}=
(P^- b P^+-P^+   a P^-)(\til a-\til b)I$.
\end{prop}
 
Hence, regarding the commutation relations for two paired operators or
transposed paired operators, we have the following.

\begin{prop}\label{prop:16jul2.6}\-\\
(i) $S_{a,b}S_{\til a,\til b}-S_{\til a,\til b}S_{a,b}=\\
\-\qquad\qquad (a-b) (P^+\til b P^--P^-\til a P^+)-
(\til a - \til b) (P^+ b P^--P^- a P^+)$;\\
(ii) $\Sigma_{a,b}\Sigma_{\til a,\til b}-\Sigma_{\til a,\til b}\Sigma_{a,b}=\\
\-\qquad\qquad (P^+ \til a P^- - P^- \til b P^+)(a-b)I-
 (P^+aP^-- P^- b P^+)(\til a - \til b) I$.
\end{prop}

The right hand sides of the equations in Propositions \ref{prop:2.2jun14} and 
\ref{prop:16jul2.6} are expressed in terms of operators of Hankel type, which allows us to
derive from these the following compactness results, where we 
 write $\GC$ for $C(\TT)$ and $\GR$ for the space of all rational functions in $L^\infty(\TT)$.
 (Standard results on Hankel operators can be found in \cite{Nik}.)

\begin{prop}\label{prop:2.4}
Suppose that $\tilde a \ne \tilde b$. Then:\\
(i) $S_{a,b}S_{\tilde a,\tilde b}-S_{a\tilde a,b\tilde b}$ is a compact
(resp. finite rank) operator if $\tilde a \in H^\infty+ \GC$, $\tilde b \in \ol{H^\infty}+\GC$
(resp.  $\tilde a \in H^\infty+ \GR$, $\tilde b \in \ol{H^\infty}+\GR$).\\
(ii) $S_{a,b}S_{\tilde a,\tilde b}-S_{\tilde a,\tilde b}S_{a,b}$ is a compact
(resp. finite rank operator) if $a,\tilde a \in H^\infty+\GC$ and $b,\tilde b \in \ol{H^\infty}+\GC$
(resp. $a,\tilde a \in H^\infty+\GR$ and $b,\tilde b \in \ol{H^\infty}+\GR$).\\
(iii) $S_{a,b}\eta I-\eta S_{a,b}$ is compact (resp. finite rank) if $\eta \in \GC$ (resp. $\eta \in \GR$).
\end{prop}

By duality, the analogous results for transposed paired operators 
are the following.

\begin{prop}\label{prop:2.5}
Suppose that $\tilde a \ne \tilde b$. Then:\\
(i) $\Sigma_{a,b}\Sigma_{\tilde a,\tilde b}-\Sigma_{a\tilde a,b\tilde b}$ is a compact
(resp. finite rank) operator if $\tilde a \in \ol{H^\infty}+ \GC$, $\tilde b \in  H^\infty+\GC$
(resp.  $\tilde a \in \ol{H^\infty}+ \GR$, $\tilde b \in  H^\infty +\GR$).\\
(ii) $\Sigma_{a,b}\Sigma_{\tilde a,\tilde b}-\Sigma_{\tilde a,\tilde b}\Sigma_{a,b}$ is a compact
(resp. finite rank operator) if $a,\tilde a \in \ol{H^\infty}+\GC$ and $b,\tilde b \in H^\infty+\GC$
(resp. $a,\tilde a \in \ol{H^\infty}+\GR$ and $b,\tilde b \in H^\infty+\GR$).\\
(iii) $\Sigma_{a,b}\eta I-\eta \Sigma_{a,b}$ is compact (resp. finite rank) if $\eta \in \GC$ (resp. $\eta \in \GR$).
\end{prop}

In particular, it is clear from (i) and (ii) in Proposition~\ref{prop:16jul2.6} that
$[S_{a,b}, S_{\til a,\til b}]=0$, that is, $S_{a,b}$ and $S_{\til a,\til b}$ commute,
if one of the following conditions hold:

\begin{eqnarray} & a=b \q &\hbox{and} \q \til a=\til b; \label{eq:16jul2.8}\\
& a,\til a \in H^\infty \q &\hbox{and} \q b,\til b \in \ol{H^\infty};\\
& a=\lambda \til a + \mu \q &\hbox{and} \q b=\lambda \til b + \mu, \q \hbox{with} \q \lambda,\mu \in \CC.\label{eq:16jul2.10}
\end{eqnarray}

The converse was provided in \cite[Thm.~4.1]{gu}, and so we have:

\begin{prop}\label{prop:16jul2.9}
$[S_{a,b},S_{\til a,\til b}]=0$ if and only if one of the conditions
\eqref{eq:16jul2.8}--\eqref{eq:16jul2.10} holds.
\end{prop}

A similar property holds for $[\Sigma_{a,b},\Sigma_{\til a,\til b}]$ by \eqref{eq:16jul2.2}.

As an illustration of the results in Proposition~\ref{prop:2.4} we have, for nondegenerate $(a,b)$, that $[S_{a,b},M_z]$ and $[S_{a,b},M_{\ol z}]$ are rank-one operators, such that
\beq\label{eq:5aug2.8}
\begin{array}{lcl}
[S_{a,b},M_z]  (f)         = & S_{a,b}zf-zS_{a,b}f                 =& (a-b)(\ol z \ol {f_-})(0), \\
{[ S_{a,b},M_{\ol z} ] (f)} =  & S_{a,b}\ol z f - \ol z S_{a,b}f    =& (b-a) \ol z f_+(0),
\end{array}
\eeq
where $f_\pm=P^\pm f$.

Analogously, $[\Sigma_{a,b},M_z]$ and $[\Sigma_{a,b},M_{\ol z}]$ are also rank-one operators,
such that
\beq\label{eq:5aug2.9}
\begin{array}{lcl}
[\Sigma_{a,b},M_z]  (f)         = & \Sigma_{a,b}zf-z\Sigma_{a,b}f                 =& [z P^-(a-b)f](\infty),\\
{[ \Sigma_{a,b},M_{\ol z} ] (f)} =  & \Sigma_{a,b}\ol z f - \ol z \Sigma_{a,b}f    =&-\ol z[P^+(a-b)f](0).
\end{array}
\eeq
Using Proposition~\ref{prop:16jul2.6} we also get the following, which is
\cite[Cor.~3.3]{CGP23}.

\begin{cor}
If $(a,b)$ is nondegenerate and $\eta \in L^\infty$ then
\[
S_{a,b}\eta I=\eta S_{a,b} \iff
\eta \in \CC \iff \Sigma_{a,b}\eta I = \eta \Sigma_{a,b}.
\]
\end{cor}

\beginpf
Note that $\eta I=S_{\eta,\eta}$ and it follows that
\[
\eta S_{a,b}=S_{a,b}\eta I \iff P^-\eta P^+= P^+ \eta P^-,
\]
which holds if and only if $\eta \in \CC$. Once again, the result for $\Sigma_{a,b}$
follows by duality.
\endpf

For $\eta  \in L^\infty\setminus \CC$ one may then ask which functions $f \in L^2$ satisfy the relation $S_{a,b}\eta f=\eta S_{a,b}f$. We have the following from
\cite[Prop.~3.4]{CGP23}.

  \begin{prop}\label{prop:2.6jun12}
  For $\eta\in L^\infty$, $(a,b)$ nondegenerate, and $f \in L^2$,
  \begin{eqnarray*}
  \eta S_{a,b}f=S_{a,b}(\eta f) & \iff& f \in \ker  H_\eta \oplus \ker \tilde H_\eta \\
  & \iff & \eta f_+ \in H^2_+, \eta f_- \in  H^2_- \\
  & \iff & \eta f_+ = P^+(\eta f) \\
  & \iff & \eta f_- = P^- (\eta f),
  \end{eqnarray*}
  where the Hankel operators $H_\eta$ and $\tilde H_\eta$ are defined by
$H_\eta=P^-\eta P^+{}_{|H^2_+}$ and $\tilde H_\eta = P^+ \eta P^-{}_{|H^2_-}$,
  and we write $f_\pm = P^\pm f$.
  \end{prop}
  
  So 
  \[
  \ker [S_{a,b},\eta I]= \ker  H_\eta \oplus \ker \tilde H_\eta.
  \]
  
  Thus, in particular, for nondegenerate $(a,b)$, the space of all functions $f \in L^2$ such that
$S_{a,b}zf=zS_{a,b}f$ is given by
\beq
\ker [S_{a,b},M_z]=\ol z H^2_- \oplus H^2_+ = \{f \in L^2: (\ol z \ol{f_-})(0)=0\}
\eeq
and, analogously,
\beq
S_{a,b} \ol z f = \ol z S_{a,b}f \iff f \in H^2_- \oplus zH^2_+ \iff f_+(0)=0.
\eeq
  
  We may prove a similar result for transposed paired kernels.
  
  \begin{prop}\label{prop:2.6Ajun12}
For $\eta\in L^\infty$, $(a,b)$ nondegenerate, and $f \in L^2$, we have
$\Sigma_{a,b}\eta f= \eta \Sigma_{a,b} f \iff (a-b)f \in \ker [S_{a,b},\eta I]$,
and hence $\ker [\Sigma_{a,b},\eta I]= 
 \{f\in L^2:(a-b)f \in \ker  H_\eta \oplus \ker \tilde H_\eta\}$.
\end{prop}
\beginpf
$P^+ a\eta f + P^- b\eta f = \eta P^+ a f + \eta P^- bf
\iff P^+ a\eta f - P^+ b\eta f = \eta P^+ a f - \eta P^+ bf$.
This holds if and only if $P^+\eta (a-b)f = \eta P^+ (a-b)f$,
and so by Proposition \ref{prop:2.6jun12} if and only if 
$(a-b)f \in \ker [S_{a,b},\eta I]$.
\endpf

Finally, one may ask when a transposed paired operator $\Sigma_{a,b}= S_{\ol a,\ol b}^*$ is
also a paired operator. We have the following answer from \cite[Prop.~2.2, Cor.~2.3]{gu},
to which we provide an alternative proof.

\begin{prop} 
$S^*_{a,b}=S_{c,d}$ for some $c,d \in L^\infty$ if and
only if $a-b \in \CC$, and, in that case, $S^*_{a,b}=S_{\ol a,\ol b}$. Thus $S_{a,b}$
is self-adjoint if and only if $a$ and $b$ are real-valued functions with $a-b \in \RR$.
\end{prop}

\beginpf
Suppose that 
\beq\label{eq:aug3A}
P^+\ol a f + P^- \ol b f = c P^+f + d P^-f \qquad \hbox{for all} \q f \in L^2.
\eeq
Then, for $f=1$, we get that $P^+\ol a + P^- \ol b =c$, so $P^+ \ol a=P^+ c$ and 
$P^-\ol b=P^- c$.
On the other hand, taking $f=\ol z$,   we get $P^+(\ol a\ol z)+ P^-(\ol b \ol z)= d \ol z    $, i.e.,
\[
\ol z (P^+ \ol a)-\ol z (P^+\ol a)(0)+\ol z(P^- \ol b)+ \ol z(P^+ \ol b)(0)=d \ol z,
\]
that is, $\ol z  c-\ol z k=\ol z d$ with $k= (P^+(\ol a-\ol b))(0) \in \CC$.

Therefore, $c=d+k$ and it follows from \eqref{eq:aug3A} that
\[
P^+(\ol a-k)I+P^-\ol b I=dI
\]
so, by Corollary~\ref{cor:14aug2.1A}, $\ol a-k=\ol b=d=c-k$; hence
$a-b \in \CC$ and $S^*_{a,b}=S_{c,d}$. The converse is straightforward.
\endpf

Turning our attention to kernels, and given the close connection between paired operators of the form $S_{a,b}$ and transposed paired operators
of the form $\Sigma_{a,b}$, which was mentioned in the introduction and the beginning
of this section, one might expect their kernels to have very similar properties. Moreover, we have the following simple relation, which generalizes \cite[Prop.~6.4]{CGP23}.

\begin{prop}\label{prop:6.4CGP23}
For nondegenerate $(a,b)$
there is a well-defined injective linear mapping  $ \GJ:\ker \sigab \to \ker S_{a,b}$
given by 
\beq  \GJ\psi=(a-b)\psi.
\eeq
If there exist $a',b' \in L^\infty$ such that
$aa'+bb'=1$, then $\GJ$ is an isomorphism with inverse
\[
\GJ^{-1} = a'P^- - b' P^+: \ker S_{a,b} \to \ker \Sigma_{a,b},
\]
and in that case 
\[
(a-b)\ker \Sigma_{a,b}=\ker S_{a,b}.
\]
\end{prop}
\beginpf
For the last part, note that $a'P^-(a-b)\psi-b'P^+(a-b)\psi= a'(a\psi)+ b'(b\psi)=\psi$
when $P^+(a\psi)=P^-(b\psi)=0$, i.e., when $\psi \in \ker \Sigma_{a,b}$.
\endpf

Somewhat surprisingly, in spite of some similarities, the two types of kernels present
significantly different and even contrasting properties, as shown in the next two sections,
where we study each type of kernel separately.

\section{Paired kernels}
\label{sec:3}

We start by noting that the kernels of paired operators, called {\em paired kernels}, are determined
by the solutions of a Riemann--Hilbert problem for the unknowns $\phi_\pm=P^\pm \phi$:
\beq\label{eq:17jul3.1}
\phi \in \ker S_{a,b} \iff a\phi_+ + b \phi_- = 0.
\eeq
It is clear from \eqref{eq:17jul3.1} that paired kernels are closely connected with Toeplitz kernels. Indeed, defining the
Toeplitz operator
\beq
T_{a/b}: \GD_{a/b}:= \{f \in H^2_+: (a/b)f \in L^2 \} \to H^2_+, \label{eq:17jul3.2}
\eeq
\beq
f \mapsto P^+ (a/b)f, \label{eq:17jul3.3}
\eeq
which is bounded on $H^2_+$ if $a/b \in L^\infty$, we have from \eqref{eq:17jul3.3} that
\beq
P^+ \ker S_{a,b} = \ker T_{a/b}.
\eeq
Nevertheless, no nonzero paired kernel can be a Toeplitz kernel, or even be contained in $H^2_+$.
This is a consequence of the following \cite[Rem.~4.2]{CP24}.

\begin{prop}\label{prop:3.5}
If $\phi \in \ker S_{a,b}$ then
$\phi=0 \iff \phi_+=0 \iff \phi_-=0$.
Therefore,
\beq
\ker S_{a,b} \cap H^2_+ = \ker S_{a,b} \cap H^2_- = \{0\}.
\eeq
\end{prop}
Clearly, for any measurable complex-valued function $\eta$ defined a.e.\ on $\TT$,
such that $\eta \ne 0$ a.e.,  we have that
\[
a\phi_+ + b \phi_- =0 \iff a\eta \phi_+ + b \eta \phi_-=0
\]
so, for any $\eta \in L^\infty$ with $\eta \ne 0$ a.e.,
\beq\label{eq:18jul3.6}
\ker S_{a,b}=\ker S_{a\eta,b\eta}.
\eeq
The next results show that this is the only case where two paired kernels are equal.

\begin{prop} \cite[Prop.~4.3]{CGP23}
If $\ker S_{a,b}  \ne \{0\}$ then, for any nondegenerate $\til a,\til b$ satisfying the  assumption \eqref{1.12} we have
\beq
\ker S_{a,b} = \ker S_{\til a,\til b} \iff a\til b=\til a b.
\eeq
\end{prop}
As a consequence, $\ker S_{a,b}$ is uniquely determined by the quotient $a/b$, which is the same
as the quotient $-\phi_-/\phi_+$ for any nonzero element of $\ker S_{a,b}$, by \eqref{eq:17jul3.1}. Thus we get the following striking property of paired kernels \cite[Thm.~4.6]{CGP23}.

\begin{thm}\label{thm:cgp4.6}
For each $\phi \in L^2 \setminus \{0\}$ there is one and only one paired kernel to which
$\phi$ belongs. If $\phi_+=I_+O_+$ and $\phi_-=I_-O_-$, where $I_+,\ol{I_-}$ are inner functions
and $O_\pm$ outer in $H^2_\pm$, then
\beq
\phi \in \ker S_{a,b}, \q \hbox{with} \q a=\ol{I_+}H_{1+}\ol{h_2+}, \q b=-z\ol{I_-}\, \ol{h_{1+}} H_{2+},
\eeq
where $H_{1+},H_{2+}, h_{1+},h_{2+} \in H^\infty$ are such that
\beq
O^{-1}_+ = \frac{H_{1+}}{H_{2+}}, \q z\ol{O^{-1}_-} = \frac{h_{1+}}{h_{2+}}.
\eeq
\end{thm}

\beginpf
To derive the formulae for $a$ and $b$ we note that
\[
\frac{\ol{I_+}}{O_+}\phi_+= \frac{\ol{I_-}}{O_-}\phi_-.
\]
Now $O_+^{-1}$ belongs to the Smirnoff class $\mathcal N_+$ so there are $H_{1+}$
and $H_{2+} \in H^\infty$ such that $O_+^{-1}  = H_{1+}/H_{2+}$; similarly
$z \ol{O_-}^{-1}= h_{1+}/h_{2+}$ with $h_{1+},h_{2+} \in H^\infty$.
The rest is clear.
\endpf

\begin{cor}
If $\ker S_{a.b} \cap \ker S_{\til a,\til b} \ne \{0\}$, then
$\ker S_{a,b}=\ker S_{\til a,\til b}$.
\end{cor}

\begin{cor}
$\ker S_{a,b} \subseteq \ker S_{\til a,\til b}
\iff \ker S_{a,b}= \{0\}$ or $\ker S_{a,b} = \ker S_{\til a,\til b}$.
\end{cor}

When do we have  $\ker S_{a,b} \ne \{0\}$? From the proof of \cite[Thm.~4.6]{CGP23}
we   obtain a necessary and sufficient condition for a nonzero kernel.

\begin{prop}\label{prop:19jul3.6}
$\ker S_{a,b} \ne \{0\}$ if and only if 
\beq\label{eq:19jul3.10}
\frac{a}{b}= I_- O_- O_+^{-1},
\eeq
where $I_+$ is inner and $O_\pm$ are outer in $H^2_\pm$.
\end{prop}

\beginpf
If \eqref{eq:19jul3.10} holds then $a O_+- b I_-O_-=0$, so that $O_+ - I_-O_ -\in \ker S_{a,b}$.

Conversely, if $a \phi_+ + b \phi_- =0$, then by taking inner--outer factorizations
we easily obtain  \eqref{eq:19jul3.10}.
\endpf

Note that this coincides with the necessary and sufficient
condition for $\ker T_{a/b}$ to be nontrivial, when $a/b \in L^\infty$, obtained in
\cite{MP05,CP18}.

An analogue of Coburn's lemma for Toeplitz kernels \cite{coburn}
also holds for paired kernels. To see this, however, it is useful to reformulate Coburn's
lemma, usually expressed as:
\[
\hbox{if} \q g \ne 0 \q \hbox{then either} \q \ker T_g = \{0\} \q \hbox{or} \q \ker T^*_g = \{0\}.
\]
Note that, while the adjoint of a Toeplitz operator is also a Toeplitz operator, the adjoint of a paired operator
is not, in general, a paired operator of the same type. Noting, however, that we
can always assume that the symbol in $\ker T_g$ satisfies $|g|=1$ \cite{hayashi85}
and, on the other hand,
\[
\ker T_g \cong \ker S_{g,1}, \q \ker T^*_g= \ker T_{\ol g} \cong \ker S_{\ol g,1}=\ker S_{1,g}
\]
by \eqref{1.11}, we can also formulate Coburn's lemma, equivalently,
as follows, in terms of paired operators.

\begin{thm} (Coburn's lemma) If $|g|=1$ then either $\ker S_{g,1}=\{0\}$
or $\ker S_{1,g}=\{0\}$.
\end{thm}

We can formulate the following analogues for paired kernels, which slightly generalize Theorem 6.8 in \cite{CGP23}.

\begin{thm}\label{thm:19jul3.8}
(i) Either $\ker S_{a,b}=\{0\}$ or $\ker S_{b,a}=\{0\}$.\\
(ii) If there exist $a',b' \in L^\infty$ such that $aa'+bb'=1$, then either
$\ker S_{a,b}=\{0\}$ or $\ker S^*_{a,b}=\{0\}$.
\end{thm}

\beginpf
(i) Suppose that there exist $\phi,\psi \ne 0$ with $\phi\in \ker S_{a,b}$ and
$\psi \in \ker S_{b,a}$. Then we have
\[
a\phi_+=-b\phi_- \qquad \hbox{and} \qquad b\psi_+ = -a\psi_-.
\]
Therefore $ab \phi_+ \psi_+ = ab \phi_- \psi_-$, and, since $a,b \ne 0$ a.e.\ on $\TT$,
it follows that
\[
\underbrace{\phi_+ \psi_+}_{\in H^1} = \underbrace{\phi_- \psi_-}_{\in \overline{H^1_0}},
\]
and so both terms above are $0$.
This implies that either $\phi_+=0$ or $\psi_+=0$ and hence, by Proposition~\ref{prop:3.5},
$\phi=0$ or $\psi=0$, which is a contradiction.\\
(ii) This follows from Proposition~\ref{prop:6.4CGP23}, since 
\[
\ker S^*_{a,b}=\ker \Sigma_{\ol a,\ol b}=\GJ^{-1} \ker S_{\ol a,\ol b}= \GJ^{-1}(J \ker S_{b,a}),
\]
where $\GJ$ was defined in Proposition~\ref{prop:6.4CGP23} and $J: \ker S_{b,a} \to \ker S_{\ol a,\ol b}$, defined by $J\phi=\ol z\ol \phi$, is a bijective antilinear operator. 
\endpf

\section{Transposed paired kernels}
\label{sec:4}

It is not difficult to see that
\begin{eqnarray}
\ker \Sigma_{a,b} &=& \{\phi \in L^2: P^+a\phi=P^-b\phi=0\} \label{eq:17jul4.1} \\
&=& \{\phi \in L^2: a\phi=\psi_-, \q b\phi=\psi_+ \q \hbox{with}\q  \psi_\pm \in H^2_\pm \}
\\
&=& L^2 \cap a^{-1}H^2_- \cap b^{-1} H^2_+, \q \hbox{for non-degenerate }   (a,b).
\label{eq:17jul4.3}
\end{eqnarray}
Thus, in contrast with paired kernels, the elements of a transposed paired kernel are not, in general, determined by
the solution of a scalar Riemann--Hilbert problem as in \eqref{eq:17jul3.1}.
As a consequence, transposed paired kernels possess properties that are different, and even
in contrast, with those of paired kernels and they have different relations with Toeplitz kernels.

Defining $T_{a/b}$ as in \eqref{eq:17jul3.2} and \eqref{eq:17jul3.3}, and taking into
account the fact that
\[
a^{-1}H^2_- \cap b^{-1} H^2_+ = b^{-1}\left( \frac{b}{a}H^2_- \cap H^2_+ \right),
\]
we have, from \eqref{eq:17jul4.3}:

\begin{prop}
(i) $\ker \Sigma_{a,b}= \left( \dfrac{1}{b} \ker T_{a/b} \right) \cap L^2$;\\
(ii) 
\beq
b \ker \Sigma_{a,b} \subseteq \ker T_{a/b}
\eeq
and, if $b \in \G L^\infty$, then
\beq\label{eq:17jul4.5}
b \ker \Sigma_{a,b}=\ker T_{a/b}.
\eeq
\end{prop}
In particular, we see from \eqref{eq:17jul4.5} that
\beq
b \in \G H^\infty \implies \ker \Sigma_{a,b} = \frac{1}{b} \ker T_{a/b} = \ker T_a,
\eeq
so, in contrast with the paired kernels studied in Section \ref{sec:3},
nontrivial transposed paired kernels can be Toeplitz kernels. Indeed, we have the following.

\begin{prop}
Every Toeplitz kernel can be written as
\beq
\ker T_a=\ker \Sigma_{a,b} \qquad \hbox{for any} \q b \in H^\infty, \q \hbox{outer}.
\eeq
\end{prop}
\beginpf
If $b \in H^\infty$ is outer, then $b^{-1}H^2_+ \cap L^2 = H^2_+$ and it
follows from \eqref{eq:17jul4.3} that $\ker \Sigma_{a,b}=a^{-1}H^2_- \cap H^2_+ = \ker T_a$.
\endpf

\begin{cor}
If $a \in \ol{H^\infty}$ and $b \in H^\infty$, with $b$ outer, then $\ker S_{a,b}$
is a model space $K_\theta$, where $\theta$ is the inner factor of $\ol a$. In
particular,
\beq
K_\theta = \ker \Sigma_{\ol\theta,1}.
\eeq
\end{cor}

These results naturally raise the question of how to characterise $\ker \Sigma_{a,b} \cap H^2_\pm$, and, in particular, when $\ker \Sigma_{a,b}$ is contained in $H^2_\pm$
or vice-versa. The answer
reveals some interesting connections of transposed paired kernels with Toeplitz and Hankel operators.

\begin{prop}\label{prop:17jul4.4}{\,}\\
(i) 
\beq \ker \Sigma_{a,b} \cap H^2_+ = \ker T_a \cap \ker H_b ,\eeq
(ii)
\[ \ker \Sigma_{a,b} \cap H^2_- = \ker \wt H_a \cap \ker \wt T_b.
\]
\end{prop}
\beginpf
(i) and (ii) follow straightforwardly from \eqref{eq:17jul4.1}.
\endpf

The following lemma will be used several times in what follows.

\begin{lem}\label{prop:17jul4.5}
If $\ker \Sigma_{a,b} \ne \{0\}$ then there exist $\phi \in \ker \Sigma_{a,b}$ such that $a\phi=\psi_-$ with $\psi_- \in H^2_-$, outer, and $\check \phi \in \ker \Sigma_{a,b}$ such that
$b \check\phi = \psi_+$ with $\psi_+ \in H^2_+$, outer.
\end{lem}

\beginpf
Let $0 \ne \til\phi \in \ker \Sigma_{a,b}$; then $a\til\phi=\psi_-$ and $b\til \phi=\til \psi_+$
where $\til\psi_\pm \in H^2_\pm$. Let $\til\psi_-=I_-\psi_-$, where $\ol{I_-}$ is inner and $\psi_-$ is outer in $H^2_-$. Then we have
$a(\ol{I_-}\til\phi)=\psi_-$, $b(\ol{I_-}\til\phi)=\ol{I_-}\til \psi_+ \in H^2_+$, so
$\phi:= \ol{I_-}\til\phi$ satisfies $a\phi=\psi_-$ with $\psi_- \in H^2_-$, outer.
We can prove similarly the second part of the proposition.
\endpf

\begin{prop}\label{prop:17jul4.6}
$\ker \Sigma_{a,b} \subseteq H^2_+ \iff \ker\Sigma_{a,b}=\ker T_a \cap \ker H_b$,
and this occurs if and only if there exists $\psi^o_+ \in H^2_+$, outer, such that
$\psi^o_+/b \in H^2_+$.
\end{prop}

\beginpf
The first equivalence is a consequence of Proposition~\ref{prop:17jul4.4}~(i).
For the second equivalence, suppose that $\psi_+^o/b \in H^2_+$ with $\psi^o_+ \in H^2_+$,
outer. Then for any $\phi \in \ker \Sigma_{a,b}$ we have that $b\phi=\psi_+ \in H^2_+$,
so 
\[
\phi = \frac{\psi_+}{b}= \frac{\psi_+ (\psi_+^o/b)}{\psi_+^o} \in \GN^+ \cap L^2 = H^2_+.
\]
(Here $\GN_+$ is the Smirnoff class \cite[Thm.~2.11]{duren}.)
\\
Conversely, if $\ker \Sigma_{a,b} \subseteq H^2_+$, then by Lemma~\ref{prop:17jul4.5}
there exists $\psi^o_+ \in H^2_+$, outer, such that $b \check\phi = \psi^o_+$ and
$\check\phi \in \ker \Sigma_{a,b} \subseteq H^2_+$,
so $\psi^o_+/b \in H^2_+$.
\endpf

\begin{prop}
\beq\label{eq:17jul4.10}
H^2_+ \subseteq \ker\Sigma_{a,b} \iff a=0 \q \hbox{and} \q b \in H^\infty
\eeq and
\[
H^2_+ = \ker \Sigma_{a,b} \iff a=0,\  b \in H^\infty \q \hbox{with}\q  b   \hbox{ outer}.
\]
\end{prop}

\beginpf
If $H^2_+ \subseteq \ker \Sigma_{a,b}$ then
$H^2_+ = \ker \Sigma_{a,b} \cap H^2_+ = \ker T_a \cap \ker H_b$.
It follows that $\ker T_a=H^2_+$, so $a=0$ and $\ker H_b= H^2_+$, so
$b \in H^\infty$. The converse is obvious. The conditions for equality follow from 
\eqref{eq:17jul4.10} and Proposition~\ref{prop:17jul4.6}.
\endpf

Analogous results hold if we replace $H^2_+$ with $H^2_-$, showing that
we can also have $\ker \Sigma_{a,b} \subseteq H^2_-$,
$H^2_- \subseteq \ker \Sigma_{a,b}$ and $\ker \Sigma_{a,b}=H^2_-$.
In particular, if $h_+ \in H^\infty$ is outer, then
\beq
\ker \Sigma_{\ol{h_+},\theta}=\ol z \ol{K_\theta} = \ol\theta K_\theta \subseteq H^2_-.
\eeq
Now,  it is clear that, in general, for $\eta \in L^\infty$,
\beq
\ker \Sigma_{a,b} \ne \ker \Sigma_{a\eta,b\eta},
\eeq
in contrast with \eqref{eq:18jul3.6}. The question then is when two transposed kernels are equal, or related by inclusion.
Here, again, the results are in contrast with those previously obtained
for paired kernels. First, however, we prove the following.

\begin{prop}
If $\phi \in \ker \Sigma_{a,b}$ and $\phi \ne 0$, then $\phi$ cannot vanish on a subset of $\TT$
with positive measure.
\end{prop}

\beginpf
If $\phi$ vanishes on a set $E \subset \TT$ with poisitve measure, then
$a \phi=\psi_-=0$ (since it is in $H^2_-$ and vanishes on $E$), and $b\phi=\psi_+=0$
(similarly), so $\phi=0$.
\endpf

\begin{prop}
Let $\ker \Sigma_{a,b} \ne \{0\}$. Then:
\\
(i) $\ker \Sigma_{a,b} \subseteq \ker \Sigma_{\til a,\til b}$ if and only if
$\til a/a=\til\psi_-/\psi_-$ with $\psi_-,\til \psi_- \in H^2_-$ and $\psi_-$ outer, and
$\til b/b = \til \psi_+/\psi_+$ with $\psi_+,\til \psi_+ \in H^2_+$ and $\psi_+$ outer;\\
(ii) $\ker \Sigma_{a,b}=\ker \Sigma_{\til a,\til b}$ if and only if $\til a/a$ and $\til b/b$ can be
written as in (i) with $\til\psi_-$ and $\til \psi_+$ outer in $H^2_-$ and $H^2_+$, respectively;\\
(iii) if $h_- \in \ol{H^\infty}$ and $h_+ \in  H^\infty$, then
\[
\ker \Sigma_{a,b} \subseteq \ker \Sigma_{ah_-,bh_+},
\]
where the inclusion is strict if $\ol{h_-}$ or $h_+$
has a nonconstant inner factor, and  equality if $h_-$ and $h_+$ are outer functions.
\end{prop}

\beginpf
(i) Let $\{0\} \ne \ker\Sigma_{a,b} \subseteq \ker \Sigma_{\til a,\til b}$ and let
$\phi\in \ker \Sigma_{a,b}$ with $a\phi = \psi_- \in H^2_-$, $\psi_-$ outer
(see Lemma~\ref{prop:17jul4.5}). Then we have
\[
a\phi= \psi_-, \q \til a \phi= \til\psi_- \in H^2_-,
\]
so $\til a/a=\til \psi_-/\psi_-$ with $\psi_-,\til\psi_- \in H^2_-$ and $\psi_-$ outer.
Analogously, using Lemma~\ref{prop:17jul4.5}, we can write $\til b/b=\til \psi_+/\psi_+$ with
$\psi_+,\til \psi_+ \in H^2_+$ and $\psi_+$ outer.

Conversely, if $\til a/a=\til\psi_-/\psi_-$ and $\til b/b = \til\psi_+/\psi_+$, where
$\til \psi_\pm,\psi_\pm \in H^2_\pm$, with $\psi_\pm$, and $\phi\in \ker\Sigma_{a,b}$, then 
$a\phi=\phi_- \in H^2_-$, $b\phi=\phi_+ \in H^2_+$, and
\begin{eqnarray*}
\til a \phi&=& \frac{\til a}{a}(a\phi)= \frac{\til \psi_-}{\psi_i} \in \ol z \ol{\GN_+} \cap L^2 = H^2_-,\\
\til b \phi &=& \frac{\til b}{b} (b\phi) = \frac{\til \psi_+}{\psi_+} \phi_+ \in \GN_+ \cap L^2 = H^2_+
\end{eqnarray*}
where $\GN_+$ is the Smirnoff class \cite[Thm.~2.11]{duren}. So $\phi \in \ker \Sigma_{\til a,\til b}$.

(ii) If $\psi_-,\til\psi_- \in H^2_-$ are outer and $\psi_+,\til\psi_+ \in H^2_+$ are outer, then by (i) we have
$\ker\Sigma_{a,b} \subseteq \ker \Sigma_{\til a,\til b}$ and $\ker \Sigma_{\til a,\til b} \subseteq \ker \Sigma_{a,b}$.

Conversely, assume that $\til a/a=\til\psi_-/\psi_-$ and $\til b/b=\til\psi_+/\psi_+$ as in (i), where, for instance $\til\psi_- \in H^2_-$ is not outer, so
$\til a/a=\ol\theta \eta_-/\psi_-$ with $\theta$ inner and $\eta_-,\psi_- \in H^2_-$ outer.
Then by (i) we have that $\ker\Sigma_{a,b} \subseteq \ker\Sigma_{\til a,\til b}$, but we cannot have
$\ker\Sigma_{\til a,\til b} \subseteq \ker\Sigma_{a,b}$.
Indeed, let $\psi \in \ker \Sigma_{\til a,\til b}$ be such that (see Lemma~\ref{prop:17jul4.5})
$\til a \phi=\phi_-$ with $\phi_-$ outer in $H^2_-$. 
If $\phi\in \ker\Sigma_{a,b}$ then we also have $a\phi=\chi_- \in H^2_-$; from
this it follows that
\[
\til a\phi=a \ol\theta \frac{\eta_-}{\psi_-} \phi=\phi_- \q \hbox{and} \q a\phi=\chi_-,
\]
so that $\phi_-= \ol \theta\dfrac{\eta_-}{\psi_-}\chi_-$, which is impossible because $\phi_-$ is outer in $H^2_-$.

(iii) This is a consequence of (i) and (ii).
\endpf

When is $\ker\Sigma_{a,b} \ne\{0\}$? We have the following necessary and sufficient condition.

\begin{prop}
$\ker \Sigma_{a,b} \ne \{0\}$ if and only if there exist $I_-,O_-,O_+$ such that $\ol{I_-}$ is inner, $O_\pm \in H^2_\pm$
are outer and
\beq\label{eq:18jul4.13}
\frac{a}{b}= I_- O_- O^{-1}_+, \q \frac{O_-}{a}, \frac{O_+}{b} \in L^2.
\eeq
\end{prop}

\beginpf
Let $\phi\in \ker\Sigma_{a,b}$, $\phi\ne 0$. Since $a\phi=\psi_- \in H^2_-$ and $b\phi=\psi_+ \in H^2_+$,
we have that $a/b=\psi_-/\psi_+$ can be written as in \eqref{eq:18jul4.13}, where $\psi_-=J_-O_-$ and
$\psi_+=J_+O_+$ with $\ol{J_-},J_+$ inner and $O_\pm \in H^2_\pm$ outer, and then $I_-=J_-/J_+$.

We then have $O_-/a=\ol{J_-}\phi \in L^2$ and $O_+/b= \ol{J_+}\phi \in L^2$.

Conversely, if \eqref{eq:18jul4.13} holds, then
\begin{eqnarray*}
a \frac{O_+}{b} &=& (I_- O_- O^{-1}_+ b) \frac{O_+}{b} = I_- O_- \in H^2_- ,\\
b \frac{O_+}{b} &=& O_+ \in H^2_+,
\end{eqnarray*}
so $0 \ne O_+/b \in \ker \Sigma_{a,b}$.
\endpf

\begin{rem}{\rm
Comparing \eqref{eq:18jul4.13} with \eqref{eq:19jul3.10} in Proposition~\ref{prop:19jul3.6},
and taking Proposition~\ref{prop:6.4CGP23} into account, we see that if $\ker \Sigma_{a,b} \ne \{0\}$ then $\ker S_{a,b} \ne \{0\}$. However, we may have
$\ker \Sigma_{a,b}=\{0\}$ and $\ker S_{a,b} \ne \{0\}$.
For example,
let $a(z)=1+1/z$ and $b(z)=z+1$. Then with $f(z)=1-1/z$ we have $f \in \ker S_{a,b}$ since
\[
(1+1/z)P^+f + (z+1) P^- f = 1+1/z - (z+1)/z=0.
\]
But $\ker \Sigma_{a,b} = \{0\}$, since
if $P_+(af)=P_-(bf)=0$ then $(1+1/z)f  \in H^2_-$ and $ (z+1)f  \in H^2_+$,
which implies that $(z+1) f \in zH^2_- \cap H^2_+ = \CC$. But this
is impossible unless $f=0$, since $f \in L^2(\TT)$.

Indeed, in this case $a/b=\ol z$ (i.e. $I_-=\ol z$, $O_-=O_+=1$), but
$1/a, 1/b \not\in L^2$.
}
\end{rem}

Another striking difference between paired kernels and transposed paired kernels is expressed by the following property.

\begin{prop} Not every $\phi \in L^2$ belongs to a transposed paired kernel.
 In particular, if $\phi$ vanishes on a set of positive measure in $\TT$,
 then $\phi$ does not belong to $\ker \sigab$ for any $a,b$, unless $\phi=0$.
 \end{prop}
 
 \beginpf
 If $\phi \in \ker\sigab$ then $P^+a\phi=P^-b\phi=0$. Thus
 $\phi=h_+/b$ for some $h_+ \in H^2_+$ and, if $\phi$
 vanishes on a set of positive measure, then $h_+=0$.
 \endpf
 
 However we do get an analogue of Coburn's lemma in the case of transposed paired kernels, as in the case of paired kernels.

\begin{prop}
If $a,b \ne 0$ a.e., then either
$\ker \sigab$ or $\ker \Sigma_{b,a}$ is $0$.
\end{prop}

\beginpf
This follows from Proposition \ref{prop:6.4CGP23} and Theorem~\ref{thm:19jul3.8}.
\endpf

\section{Invariance properties of paired kernels}
\label{sec:5}

It is easy to see that if $X$ and $T$ are commuting operators on a Banach space $\GX$
then $\ker T$ is invariant for $X$.
Thus from the commutation relations of Section \ref{sec:2} 
(namely Proposition \ref{prop:16jul2.9})
we have the
following.

\begin{prop} \cite[Prop. 5.1]{CGP23}
If $a,\tilde a \in H^\infty$ and $b,\tilde b \in \ol{H^\infty}$,
then $\ker S_{a,b}$ is invariant for $S_{\tilde a,\tilde b}$. If
$a,\tilde a \in \ol{H^\infty}$ and $b,\tilde b \in H^\infty$,
then $\ker \Sigma_{a,b}$ is invariant for $\Sigma_{\tilde a,\tilde b}$. 
\end{prop}

In particular, since for any inner function $\theta$, we have $\ker \Sigma_{\ol\theta,1}=K_\theta$ and, on the other hand,
$\Sigma_{\ol\theta, 1}\Sigma_{\tilde a,\tilde b}=\Sigma_{\tilde a,\tilde b}\Sigma_{\ol\theta,1}$ if
$\tilde a \in \ol{H^\infty}$ and $\tilde b \in H^\infty$,
we have the following generalization of 
the well-known $S^*$-invariance property of model spaces.

\begin{cor}
Model spaces are invariant subspaces for the transposed
paired operators $\Sigma_{\tilde a,\tilde b}$ with
$\til a \in \ol{H^\infty}$ and $\til b \in H^\infty$.
\end{cor}

Together with the notion of invariant subspaces for an operator $T$, one
can define also the concept of almost-invariant 
subspaces for $T$.

\begin{defn}
Let $\GH$ be a Hilbert space, $X \in \GL(\GH)$, and let $\GM$ be a subspace of 
$\GH$. We say that $\GM$ is {\em almost invariant with defect $m$}
for $X$ if and only if
\beq
X\GM \subset \GM \oplus \GF,
\eeq
where $\GF$ is a finite-dimensional subspace of minimal dimension $m$.
\end{defn}
It was shown by Popov and Tcaciuc \cite{PT13}
that every Hilbert space operator admits an almost-invariant half-space (that is, a subspace of infinite dimension
and infinite codimension) of defect 1. 
 The following result from \cite{CKP24}
also holds.

\begin{prop}
Let $X,T \in \GL(\GH)$, where $\GH$ is a Hilbert space,  such that $XT-TX$ is a finite-rank operator with rank $r$. Then
$\ker T$ is almost invariant with defect $m \le r$ for $X$, i.e., 
\beq\label{eq:5aug5.1A}
f \in \ker T \implies Xf \in \ker T \oplus \GF,
\eeq
where $\GF$ is a subspace of dimension at most $r$.
\end{prop}

From this proposition and Propositions \ref{prop:16jul2.6}--\ref{prop:2.5} we then have the following.

\begin{prop}
If $a,\tilde a \in H^\infty+\GR$ and $b,\tilde b \in \ol{H^\infty}+\GR$, then
$\ker S_{a,b}$ is almost invariant for $S_{\tilde a,\tilde b}$ with defect less than or equal to $\rank [S_{a,b},S_{\tilde a,\tilde b}]$.
\end{prop}

A similar property holds for transposed paired operators. From \eqref{eq:5aug2.8} and 
\eqref{eq:5aug2.9} we have then:

\begin{cor}
$\ker S_{a,b}$ and $\ker \Sigma_{a,b}$ are almost invariant
for multiplication by $z$ and $\ol z$ with defect 1.
\end{cor}

Besides invariant and almost-invariant subspaces for an operator $T$, one
can also define nearly-invariant subspaces. However, while the former
are defined in a general Banach space setting, the notion of
nearly $S^*$-invariant subspaces appeared in the particular context of
$H^2_+$ in the work of Hitt, Hayashi and Sarason
\cite{hayashi,hitt,sarason},
motivated by questions of invariant subspaces and kernels of
Toeplitz operators.

\begin{defn}
Let $\GM$ be a proper subspace of $H^2_+$. Then we say that $\GM$ is  {\em nearly $S^*$-invariant\/} if and only if
\beq\label{eq:nivdef}
f \in \GM, f(0)=0 \implies S^*f \in \GM.
\eeq
\end{defn}

It is well known that kernels of Toeplitz operators   are nearly $S^*$- invariant.
However, kernels of Hankel operators $H_g: H^2_+ \to H^2_-$
with $H_g \phi = P^- g\phi$ and  
truncated Toeplitz operators $A_g^{\theta}: K_\theta \to K_\theta$
with $A_g^{\theta} \phi=P_\theta g\phi$, where $P_\theta$ is the orthogonal
projection onto $K_\theta$, are not in general
nearly $S^*$-invariant, and the following weaker notion is of interest in those cases.

\begin{defn}\label{def:5aug5.5}\cite{CGP20}
Let $\GM$ be a proper subspace of $H^2_+$. Then we say that $\GM$ is {\em nearly $S^*$-invariant
with defect $m$\/} if and only if
\beq
f \in \GM, f(0)=0 \implies S^*f \in \GM \oplus \GF,
\eeq
where $\GF$ is a finite-dimensional subspace of minimal dimension $m$.
\end{defn}
A description of nearly $S^*$-invariant subspaces of defect $m$ was given in \cite{CGP20}.
A nontrivial kernel of a Hankel operator, being of the form $\theta H^2_+$ for $\theta$ inner,
is nearly $S^*$-invariant with defect 1 if $\theta(0)=0$, since if $f \in \theta H^2_+$ with
$f(0)=0$ then $\bar z f \in \theta H^2 + \CC \bar z \theta$; while if $\theta(0) \ne 0$
it is already  nearly $S^*$-invariant.
Likewise the kernel of an asymmetric truncated Toeplitz operator is nearly $S^*$-invariant
with defect 1, as   shown in \cite{ryan20}; see also \cite{CKP24}.
Similar observations apply to kernels of certain finite-rank perturbations of Toeplitz operators
\cite{LP20,DS24}.
\\

Near $S^*$-invariance, with or without positive
defect, can be seen and interpreted from various perspectives,
leading to different generalizations of the notion
of near invariance, allowing it to
apply beyond the context of $H^2_+$. 
First, we see that \eqref{eq:nivdef} is equivalent to the condition
\beq\label{eq:5aug5.4}
f \in \GM, \ol z f \in H^2_+ \implies \ol z f \in \GM.
\eeq
So a natural generalization of Definition \ref{def:5aug5.5} is the following.

\begin{defn}\cite{CP14} Let $\GM$ be a proper subspace of $H^2_+$ and let $\eta$ be a complex-valued function defined a.e.\ on $\TT$. We
say that $\GM$ is {\em nearly $\eta$-invariant} (or {\em nearly $M_\eta$-invariant}) if and only if
\beq\label{eq:5aug5.5}
f \in \GM,\ \eta f \in H^2_+ \implies \eta f \in \GM,
\eeq
and, if $\eta \in L^\infty$, we also say that $\GM$ is {\em nearly $T_\eta$-invariant}, since
\eqref{eq:5aug5.5} is equivalent to
\[
f \in \GM,\ \eta f \in H^2_+ \implies T_\eta f \in \GM.
\]
Further, if
\beq f \in \GM,\ \eta f \in H^2_+ \implies \eta f \in \GM \oplus \GF,
\eeq
where $\GF$ is a finite-dimensional subspace of minimal dimension $m$, then
we say that $\GM$ is
{\em nearly $\eta$-invariant with defect $m$}.
\end{defn}

It was shown in \cite{CP14} that Toeplitz kernels are not only nearly $\ol z$-invariant (that is, nearly $S^*$-invariant) subspaces of $H^2_+$, but also nearly $\eta$-invariant for a wide class of functions $\eta$, including all functions in $\ol{H^\infty}$ and all rational functions
without poles in the exterior of the unit disc. This property of near $\eta$-invariance provides
important information characterizing Toeplitz kernels.
For instance, it implies that a nontrivial Toeplitz kernel cannot be of the form
$\alpha K$ with $\alpha$ inner and nonconstant and $K \subset H^2_+$, since
Toeplitz kernels are nearly $\ol\alpha$-invariant. Nor can it be of the form $(z+a)K$ with
$a \in \DD \cup \TT$ and $K \subset H^2_+$, since Toeplitz
kernels are also nearly $1/(z+a)$-invariant.

To extend this notion beyond the particular setting of $H^2_+$, the following was introduced recently in \cite{CKP24}.

\begin{defn}\label{def:3jun12}
Let $\GH$ be a Hilbert space and let $\GM$, $\GN$ be closed subspaces
of $\GH$, with $\GM \subset \GN$, such that for $X \in \GL(\GH)$ we have
\beq\label{eq:5aug5.7}
f \in \GM, \ Xf \in \GN \implies Xf \in \GM.
\eeq
Then we say that $\GM$ is {\em nearly $X$-invariant with respect to $\GN$} (for short, w.r.t.\ $\GN$). 
Also, $\GM$ is said to be {\em nearly $X$-invariant w.r.t. $\GN$ with defect $m$} if
\beq 
f \in \GM, \ Xf \in \GN \implies Xf \in \GM \oplus \GF,
\eeq
with $\GF$ as above.
\end{defn}

Note that, if $\GN$ is a closed subspace and $P_\GN$ denotes the orthogonal projection onto $\GN$, then \eqref{eq:5aug5.7} is equivalent to
\[
f \in \GM, \ Xf \in \GN \implies P_\GN Xf \in \GM,
\]
so we also say that $\GM$ is nearly $P_\GN X_{|\GN}$-invariant.

Taking, for instance $X=M_{\ol z}$, $\GN=H^2_+$, $\GH=L^2(\TT)$ and $P_\GN X=S^*$, we recover \eqref{eq:5aug5.4} from \eqref{eq:5aug5.7}.\\

There are obvious relations between the notions of almost-invariant and nearly-invariant subspaces for an operator, especially in the context of $H^2_+$.
If $Xf \in \GN$ for all $f \in \GM$, then the space $\GM$ is nearly $X$-invariant with defect $m$ w.r.t.\ $\GN$
if and only if it is almost invariant with defect $m$ for $X$.
This happens in particular if $\GN=\GH$. Moreover, 
as pointed out in \cite{CGP20}, we also have that if $\GM \subset H^2_+$ is
 almost invariant for $S^*$ with defect at most
 $m$, then it is nearly $S^*$-invariant with defect at most $m$. 
 The converse is true if (i) all elements of $\GM$ vanish at $0$, or (ii) 
  there exists $f_0 \in \GM$ with $f_0(0) \ne 0$ such that
 $S^* f_0 \in \GM + \GF$.
 
 Nevertheless, even when the notion of nearly $X$-invariant subspace of $\GH$, with or without defect, coincides with that of
 almost-invariant subspace with a certain defect for the operator $X$, the two concepts
 express different points of view and address different questions on how the operator acts on a certain space.
 
 For instance, in the context of $H^2_+$, saying that a subspace $\GM$
 is almost invariant with defect $m$ for an operator $X$ can be seen as answering the question
 of where $X(\GM)$ lies, if it is not
 contained in $\GM$. Saying that $\GM$ is nearly $X$-invariant, with or without defect, can be seen as
 answering the question of which elements of $\GM$ are mapped into $\GM$
 (or $\GM \oplus \GF$) by $X$.
 Thus, considering the case of Toeplitz kernels, which are
 $S^*$-invariant if and only if the symbol is in $\ol{H^\infty}$, we may ask which
 elements $f \in \ker T_g$ are
 mapped into $\ker T_g$ by multiplication by $\ol z$, or by $S^*$, when $g \not\in \ol{H^\infty}$.
 The answer is given by $\ol z f \in H^2_+$, i.e.,
 \beq
 f \in \ker T_g, \ \ol z f \in H^2_+ \implies \ol z f \in \ker T_g,
 \eeq
 which is precisely the same as saying that $\ker T_g$ is nearly $S^*$-invariant.

 If we now consider instead the problem of how ``far'' $S^*(\ker T_g)$ is from
 $\ker T_g$ when $g \not\in \ol{H^\infty}$, in the context of $H^2_+$, or the same
 question for $\ol z \ker T_g$, in the context of $L^2$, the answer is given by saying that
 $\ker T_g$ is almost invariant with defect 1 for both
 $S^*$ (in $H^2_+)$ and $M_{\ol z}$ (in $L^2$). Thus in the latter
 case one gets information on the image of $\ker T_g$ by $S^*$, while
 in the case of near invariance one gets information describing $\ker T_g$
 itself; in particular, there must be an element of $\ker T_g$
 which does not vanish at $0$, and so $\ker T_g \not\subseteq zH^2_+$.
 
 In the more general context of $L^2$, to extend the
 previous results to paired operators, one can look at near invariance in
 the context of $\GN$-stability, or stability in $\GN$, introduced in \cite{CKP24}.
 
 \begin{defn}
 Let $\GA \ne \{0\}$ be a closed subspace of $\GN \subset \GH$ and let $X \in \GL(\GH)$. We say that $\GA$ is $\GN$-stable for $X$, or stable in $\GN$ for $X$, if $X(\GA) \subset \GN$.
 \end{defn}
 
 Thus, according to this definition, saying that $\GM$ is nearly $S^*$-invariant is equivalent to saying
 that $\GM \cap zH^2_+$ is $\GM$-stable for the operator $S^*$; $\GM$ is invariant
 for $S^*$ if and only if it is $\GM$-stable for $S^*$.
 
 When $\GM$ is the kernel of an operator $T$ on $\GN \subset \GH$,
 this notion can be related with the commutation
 relations between two operators as follows.
 It is clear that if $X$ and $T$ commute then all elements of $\ker T$ are mapped back into $\ker T$ by the operator $X$; i.e., $\ker T$ is invariant for $X$. But
 if $X$ and $T$ do not commute on their domain, then
those elements of $\ker T$ for which $X$ and $T$ do commute, i.e.,
 \[
 \GM = \{f \in \ker T: XTf=TXf\}
 \]
 are also mapped back into $\ker T$ by $X$, and so $\GM$ is $\ker T$-stable for $X$.
 
 Comparing this with \eqref{eq:nivdef}, note that
 for $g \not\in \ol{H^\infty}$,  taking $T=T_g$ and $X=S^*$,
 the elements of $\GM$, as defined above,
 are precisely those with $f(0)=0$, so one must
 indeed have that $S^*f \in \ker T_g$.
 Noting moreover that $\ol z f \in H^2_+$ means that $f \in \ker H_{\ol z}$, we see that
 in this case
 \[
 \GM= \ker T_g \cap \ker H_{\ol z}
 \]
 and near $S^*$-invariance can be expressed as saying that $\ker T_g \cap \ker H_{\ol z}$ is stable in $\ker T_g$
 for the operator $S^*$. One could express near $\eta$-invariance analogously by saying that
 $ \ker T_g \cap \ker H_\eta$ is stable in $\ker T_g$ for $M_\eta$.
 
 Applying these ideas now to the setting of paired operators in $L^2(\TT)$, we have from Propositions \ref{prop:2.6jun12} and \ref{prop:2.6Ajun12}
 that, for $\eta\in L^\infty$ and $f \in L^2$,
 \[
 \begin{array}{rclcl}
 \eta S_{a,b}f &=& S_{a,b}\eta f & \iff & f \in \ker H_\eta \oplus \ker \wt{H}_\eta\\
 \eta\Sigma_{a,b}f&=& \Sigma_{a,b}\eta f & \iff & (a-b)f \in \ker H_\eta \oplus \ker \wt H_\eta .
 \end{array}
 \]
 Thus we have the following analogue of near $\eta$-invariance.
 
 \begin{prop}
 (i) $\ker S_{a,b} \cap (\ker H_\eta \oplus \ker \wt H_\eta)$ is stable
 in $\ker S_{a,b}$ for multiplication by $\eta$; i.e.,
 \[
 f \in \ker S_{a,b} ,\ H_\eta f_+=0, \ \wt H_\eta f_-=0 \implies \eta f \in \ker S_{a,b}.
 \]
 (ii) If $f \in \ker \Sigma_{a,b}$ and $(a-b)f \in \ker H_\eta \oplus \ker \wt H_\eta$,
 then $\eta f \in \ker \Sigma_{a,b}$.
 \end{prop}

 \begin{cor}
(i) $f \in \ker S_{a,b}, \ f_+ \in \ker H_{\ol z} \implies \ol z f \in \ker S_{a,b}$,
so there exists in every non-zero paired kernel an element $f$ with $f_+(0) \ne 0$.\\
(ii) $f \in \ker S_{a,b}, \ f_- \in \ker \tilde H_z \implies zf \in \ker S_{a,b}$,
so there exists in every non-zero paired kernel an element $f$ with $z f_-(\infty) \ne 0$.
\end{cor}
 
 Note that (i) in the above corollary generalizes a well-known property of nonzero nearly $S^*$-invariant subspaces of $H^2_+$.

\end{document}